\PassOptionsToPackage{unicode}{hyperref}
\PassOptionsToPackage{hyphens}{url}
\PassOptionsToPackage{dvipsnames,svgnames,x11names}{xcolor}
\documentclass[
  10pt,
  a4paper,
  twocolumn]{article}
\usepackage{amsmath,amssymb}
\usepackage{lmodern}
\usepackage{iftex}
\ifPDFTeX
  \usepackage[T1]{fontenc}
  \usepackage[utf8]{inputenc}
  \usepackage{textcomp} 
\else 
  \usepackage{unicode-math}
  \defaultfontfeatures{Scale=MatchLowercase}
  \defaultfontfeatures[\rmfamily]{Ligatures=TeX,Scale=1}
\fi
\IfFileExists{upquote.sty}{\usepackage{upquote}}{}
\IfFileExists{microtype.sty}{
  \usepackage[]{microtype}
  \UseMicrotypeSet[protrusion]{basicmath} 
}{}
\makeatletter
\@ifundefined{KOMAClassName}{
  \IfFileExists{parskip.sty}{%
    \usepackage{parskip}
  }{
    \setlength{\parindent}{0pt}
    \setlength{\parskip}{6pt plus 2pt minus 1pt}}
}{
  \KOMAoptions{parskip=half}}
\makeatother
\usepackage{xcolor}
\usepackage[left=2cm,right=2cm,top=2.5cm,bottom=2.5cm]{geometry}
\usepackage{longtable,booktabs,array}
\usepackage{calc} 
\usepackage{etoolbox}
\makeatletter
\patchcmd\longtable{\par}{\if@noskipsec\mbox{}\fi\par}{}{}
\makeatother
\IfFileExists{footnotehyper.sty}{\usepackage{footnotehyper}}{\usepackage{footnote}}
\makesavenoteenv{longtable}
\usepackage{graphicx}
\makeatletter
\def\maxwidth{\ifdim\Gin@nat@width>\linewidth\linewidth\else\Gin@nat@width\fi}
\def\maxheight{\ifdim\Gin@nat@height>\textheight\textheight\else\Gin@nat@height\fi}
\makeatother
\setkeys{Gin}{width=\maxwidth,height=\maxheight,keepaspectratio}
\makeatletter
\def\fps@figure{htbp}
\makeatother
\newlength{\cslhangindent}
\newlength{\csllabelwidth}
\newlength{\cslentryspacingunit} 
\newenvironment{CSLReferences}[2] 
 {
  \setlength{\parindent}{0pt}
  \ifodd #1
  \let\oldpar\par
  \def\par{\hangindent=\cslhangindent\oldpar}
  \fi
  \setlength{\parskip}{#2\cslentryspacingunit}
 }%
 {}
\usepackage{calc}

\usepackage{amsmath}
\usepackage{amssymb}
\usepackage{bm}
\usepackage{mathtools}
\usepackage{interval}
\usepackage{xurl}
\usepackage[nottoc]{tocbibind}
\usepackage{subfig}
\ifLuaTeX
  \usepackage{selnolig}  
\fi
\IfFileExists{bookmark.sty}{\usepackage{bookmark}}{\usepackage{hyperref}}
\IfFileExists{xurl.sty}{\usepackage{xurl}}{} 
\hypersetup{
  pdftitle={An Approach to Ordering Objectives and Pareto Efficient Solutions},
  pdfauthor={Sebastian Hönel; Welf Löwe},
  colorlinks=true,
  linkcolor={blue},
  filecolor={Maroon},
  citecolor={Blue},
  urlcolor={blue},
  pdfcreator={LaTeX via pandoc}}

\title{An Approach to Ordering Objectives and Pareto Efficient Solutions}
\usepackage{etoolbox}
\makeatletter
\providecommand{\subtitle}[1]{
  \apptocmd{\@title}{\par {\large #1 \par}}{}{}
}
\makeatother
\subtitle{A low- or no-preference demonstration using the Viennet function}
\author{Sebastian Hönel \and Welf Löwe}
\date{May 27, 2022}

\begin{document}
\maketitle
\begin{abstract}
  Solutions to multi-objective optimization problems can generally not be compared or ordered, due to the lack of orderability of the single objectives. Furthermore, decision-makers are often made to believe that scaled objectives can be compared. This is a fallacy, as the space of solutions is in practice inhomogeneous without linear trade-offs. We present a method that uses the probability integral transform in order to map the objectives of a problem into scores that all share the same range. In the score space, we can learn which trade-offs are actually possible and develop methods for mapping the desired trade-off back into the preference space. Our results demonstrate that Pareto efficient solutions can be ordered using a low- or no-preference aggregation of the single objectives. When using scores instead of raw objectives during optimization, the process allows for obtaining trade-offs significantly closer to the expressed preference. Using a non-linear mapping for transforming a desired solution in the score space to the required preference for optimization improves this even more drastically.
\end{abstract}

{
\hypersetup{linkcolor=}
\setcounter{tocdepth}{6}
}
\intervalconfig{soft open fences}

\newcommand*\mean[1]{\overline{#1}}
\newcommand{\abs}[1]{\left\lvert\,#1\,\right\rvert}
\newcommand{\norm}[1]{\left\lVert\,#1\,\right\rVert}
\newcommand{\infdiv}[2]{#1\;\|\;#2}
\newcommand\argmax[1]{\underset{#1}{arg\,max}}
\newcommand\argmin[1]{\underset{#1}{arg\,min}}

\hypertarget{introduction}{%
\section{Introduction}\label{introduction}}

Solving multi-objective optimization problems (MOOPs, \eqref{eq:moop}) are problems of two or more conflicting problems (or objectives).

\begin{align}
  \min_{\mathbf{x}\in\mathcal{D}}\,\left\{f_1(\mathbf{x}),\dots,f_k(\mathbf{x})\right\},\text{ where $k\geq2$.}\label{eq:moop}
\end{align}

The goal is to find some \(\mathbf{x}^\star\) in the \emph{decision space} \(\mathcal{D}\) that minimizes the loss across all single objectives simultaneously.
Their solving may lead to a (possibly infinitely) large solution space \(\mathcal{S}\).
Each vector \(\mathbf{s}_i\) in this space holds a solution for each of the optimized objectives.
The set of Pareto \emph{efficient} solutions comprises those solutions that cannot be further improved.
Because vectors in the Pareto efficient solution space cannot be ordered completely (\protect\hyperlink{ref-miettinen2008}{Miettinen 2008}), they are traditionally regarded as equally desirable (in the mathematical sense, that is).
A decision maker (DM) is then consulted to pick a preferred solution.
Even worse, it may not be trivial to pick some most preferred solution, as the role of the chosen weights is misleading (\protect\hyperlink{ref-roy1996theoretical}{Roy and Mousseau 1996}).
The mathematical sense is, however, not congruent with interpreting objectives as scores. While in some cases objectives can be scaled into some uniform, dimensionless scale, the distribution of their associated losses remains unknown and \textbf{must} be assumed to be non-uniform. In other words, not each attainable loss is equally likely, which drastically exacerbates the DMs situation of having to pick \emph{and} actually obtain the desired trade-off.

We demonstrate an approach to approximating the order of Pareto efficient solutions, based on the marginal distributions of each objective.
This approach transforms each objective into a score with standard uniform distribution, which allows comparing solution vectors.
Used as either an a priori or a posteriori method, we can show that obtaining efficient solutions significantly closer to the desired trade-off is possible, using two distinct methods.
The first method significantly improves precision by replacing the raw objectives with their scores during optimization.
The second method improves the precision even more drastically, but requires the computation of some Pareto efficient solutions, in order to learn the non-linear mapping between preferences and solutions in the score space.
None of the methods require the objectives to be scaled or having to know or approximate nadir-, ideal-, or utopian-vectors.
Furthermore, by mapping objectives to scores, it allows a DM to gain insights into the density and homegenuity of the solution space, as well as to understand which trade-offs are actually possible.
We use the Viennet function (\protect\hyperlink{ref-viennet1996}{Viennet, Fonteix, and Marc 1996}) to empirically gather some results.

\hypertarget{transformation-to-scores}{%
\section{Transformation to scores}\label{transformation-to-scores}}

Given some objective \(f:\mathbb{R}^m\to\mathbb{R}\) that is subjected to minimization, we approximate its empirical distribution by uniformly drawing vectors \(\mathbf{x}\in\mathbb{R}^m\) from the decision space \(\mathcal{D}\).
We therefore treat the outcome of each objective as a random variable that follows some distribution \(\mathcal{X}\sim\bm{\beta}\), where the parameters \(\bm{\beta}\) are unknown.
Any random variable can be transformed into another random variable with standard uniform (or arbitrary other) distribution, using the \emph{probability integral transform}.
For our purposes, the standard uniform distribution suffices, as its range is \(\interval{0}{1}\), which corresponds to what one would expect from a score.
The cumulative distribution \(\operatorname{CDF}\) of a random variable expresses the probability to find a value less than or equal to \(x\).
The results of the objective are ordered, and so is its corresponding \(\operatorname{CDF}\). However, the meaning of loss and score is still reversed, in that a low loss corresponds to a low chance of observing it.
We therefore define the score for \(f\) as in \eqref{eq:score-for-f},

\begin{align}
  S_f:\mathbb{R}^m\mapsto\interval{0}{1}=1-\operatorname{CDF}_f(f(\mathbf{x}))\label{eq:score-for-f},
\end{align}

where the operation \(\operatorname{CCDF}=1-\operatorname{CDF}\) is the \emph{complementary} \(\operatorname{CDF}\).
Now, a low loss of \(f\) will yield a high score for \(S_f\), and vice versa.
Throughout this work we use empirical \(\operatorname{CDF}\)s approximated with high precision (using \(10^6\) random decision vectors), in order to map from the loss- into the score-space. This is done to achieve high numerical precision in the context of this work.
As we will show, in practice it suffices to draw considerably fewer samples (e.g., \(10^3\)) to obtain a sufficient approximation.

\hypertarget{ordering-of-solutions}{%
\subsection{Ordering of solutions}\label{ordering-of-solutions}}

Orderability of Pareto efficient solutions is something that has eluded multi-objective optimization. In practice, the DM was only left with the mathematical equal desirability.
Using \(\operatorname{CCDF}\)s, obtained solutions can be mapped into the score space \(\mathcal{S}\). Using some final scheme that expresses the weight or importance of each score for some absolute best solution, the obtained scores for each solution can be aggregated, and the space be ordered according to that.
The result of this may or may not be a unique solution.

\hypertarget{obtaining-desired-trade-offs}{%
\subsection{Obtaining desired trade-offs}\label{obtaining-desired-trade-offs}}

While preference may be expressed, it is not necessarily adhered to, neither by the solution algorithm, nor by the problem itself.
When specifying preference, a DM might obtain a solution close to that preference. However, that solution is in the objective space \(\mathcal{O}\), \textbf{not} in the score space \(\mathcal{S}\).
The latter is perhaps the only space workable for human DMs. This means that so far, expressing preference with non-uniformly distributed objectives has never led to the desired solution.

In order to obtain a trade-off close to the desired trade-off, we identify two methods.
In the first method, the solution algorithm would solve a version of the original problem, but based on the scores \eqref{eq:moop-scores}.
Since we have a one-to-one association between preference and score, this method ought to converge towards the desired trade-off more accurately.
This method usually requires a somewhat better empirical approximation of the marginal cumulative densities, or some well-fitting parametric probability distribution.
In the former case, measures for introducing smoothness into the empirical \(\operatorname{CDF}\)s have to be taken as otherwise, the gradient of the problem will be zero, making this method only applicable in gradient-free scenarios.
This method also enables the simple weighting method and the method of \emph{weighted metrics} (or compromise programming) (\protect\hyperlink{ref-zeleny1973compromise}{Zeleny 1973}).
The weighting methods require their objectives to have the same range, which we will get by using \(\operatorname{CDF}\)s.
This is an important point, as by sufficiently approximating those, we do not have to approximate the nadir-, ideal-, or utopian-vectors, which are otherwise required for proper scaling.
Especially the nadir-vector is difficult to obtain in practice (\protect\hyperlink{ref-miettinen2008}{Miettinen 2008}), so not requiring precise approximations of any of those vectors means we can dispense with the associated effort, and may also be able to use previously incompatible solution algorithms.

\begin{align}
  \max_{\mathbf{x}\in\mathcal{D}}\,\left\{S_1(\mathbf{x}),\dots,S_k(\mathbf{x})\right\}.\label{eq:moop-scores}
\end{align}

In the second method, we attempt to ``rectify'' the preference as expressed by the DM. As we have previously established, there does not exist a linear relation between some preference and its solution in the score space.
We can however approximate that non-linear relationship, too.
By uniform randomly sampling from the \emph{preference space} \(\mathcal{P}\), we first compute a sufficiently large subset of the Pareto optimal set in the score space.
Typically, this set is considerably smaller than the full set \eqref{eq:pref-bijections}.
These solutions are computed to establish a bijection \(\mathcal{P}\to\mathcal{S}\) between the spaces for the desired- and obtained trade-off (i.e., which preference leads to which trade-off in the score space).
Then, this relation is reversed and some model \eqref{eq:pref-nl-model} is fit that minimizes the deviation.

\begin{align}
  \mathcal{D}_k&=\left\{\left(\bm{\mathcal{P}}_i,\bm{\mathcal{S}}_i\right)\right\}_{i=1}^k\text{, dataset of bijections,}\label{eq:pref-bijections}
  \\[1ex]
  \mathsf{M}&=\min\,\sum_{j=i}^{k}\,(\mathbf{S}_i-\mathbf{P}_i)^2\text{, non-linear model.}\label{eq:pref-nl-model}
\end{align}

This can be done as we now have a one-to-one correspondence between preferences and solutions in the score space.
The DM is now enabled to express their preference, and the learned model corrects this preference in order to converge to an efficient solution close to the desired trade-off.

\hypertarget{the-viennet-function}{%
\section{The Viennet function}\label{the-viennet-function}}

The Viennet function consists of three objectives, each of which takes the same two parameters \(\mathbf{x}=\{x_1,x_2\}\).
It is defined as in \eqref{eq:viennet} (\protect\hyperlink{ref-viennet1996}{Viennet, Fonteix, and Marc 1996}).
The box bounds for the decision space \(\mathcal{D}\) for this problem are sometimes limited to \(-3\leq x_1,x_2\leq3\), but there is no practical difference using the slightly larger bounds as we do.

\begin{align}
  \min_{\mathbf{x}\in\mathbb{R}^2}\,&\begin{cases}
    f_1(\mathbf{x})&=0.5\left(x_1^2+x_2^2\right) + \sin{\left(x_1^2+x_2^2\right)},
    \\[1em]
    f_2(\mathbf{x})&=\frac{(3x_1-2x_2+4)^2}{8}+\frac{(x_1-x_2+1)^2}{27}+15,
    \\[1ex]
    f_3(\mathbf{x})&=\frac{1}{x_1^2+x_2^2+1}-1.1\exp{\left(-x_1^2-x_2^2\right)},
  \end{cases}\label{eq:viennet}
  \\[1ex]
  &\text{subject to }-4\leq x_1,x_2\leq4\nonumber.
\end{align}

\begin{figure}

{\centering \includegraphics[width=0.9\linewidth]{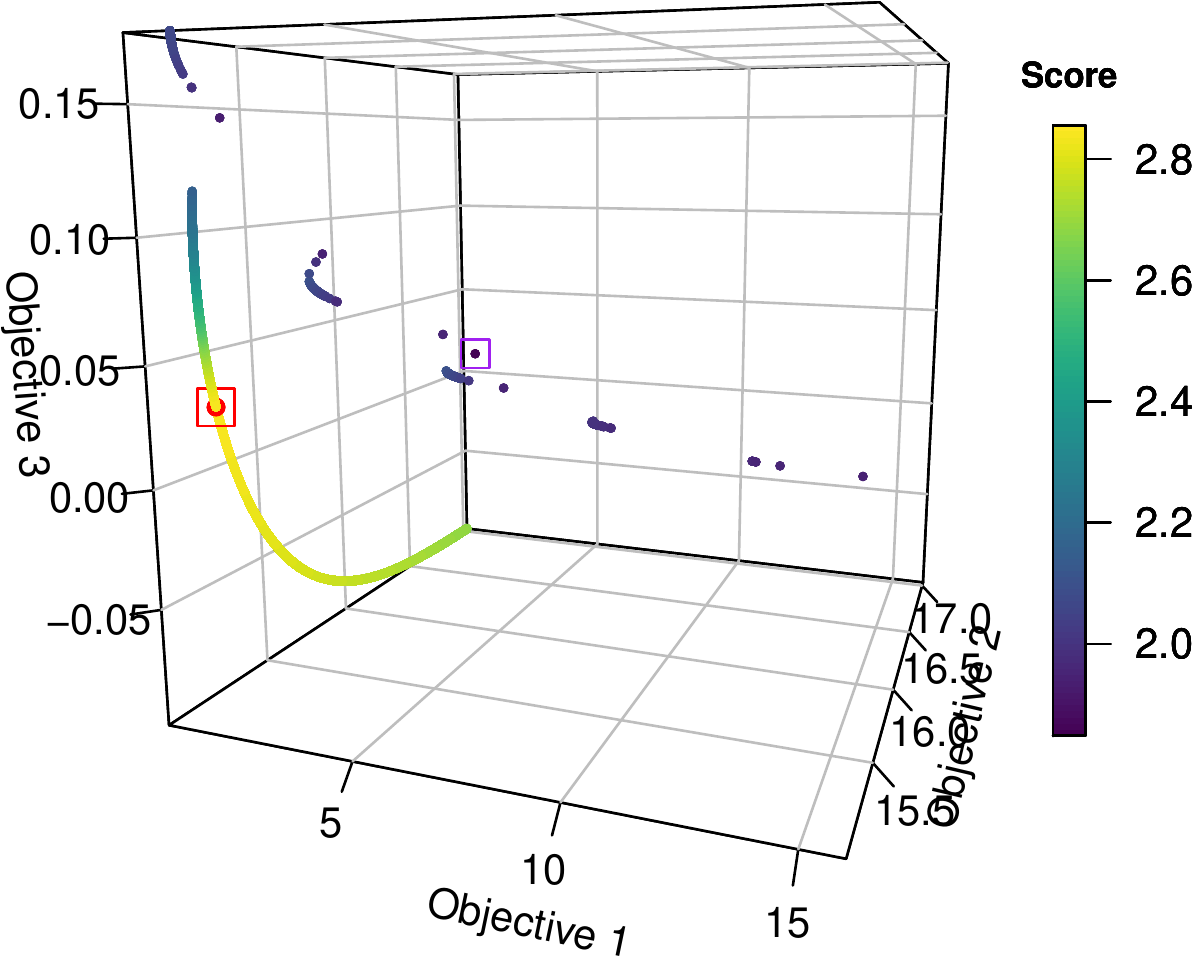} 

}

\caption{The 3-dimensional Pareto front of the Viennet function's objectives, in the respective objective spaces. The color corresponds to the sum of scores in the score space.}\label{fig:ex2-pareto-front}
\end{figure}

Figure \ref{fig:ex2-pareto-front} shows the Pareto front of the efficient set. In order to obtain the front, we compute \(50,000\) solutions by drawing preference constellations uniform randomly from the preference space \(\mathcal{P}\).
It is obvious that the three objectives do not share a common domain.
The starting point for each optimization is chosen deterministic randomly, within the box bounds of the decision space \(\mathcal{D}\).
As solution algorithm we use the derivative-free global optimization called ``DIRECT-L''.
It is based on a systematic subdivision of the search domain into increasingly smaller becoming hyperrectangles (\protect\hyperlink{ref-gablonskyK2001}{Gablonsky and Kelley 2001}).
The scalarizer used is the simple weighting method. While it has deficiencies, it can find Pareto efficient solutions if all weights are \(>0\). The weighting method requires all objectives to share some common uniform dimensionless scale, which is fulfilled when using scores instead of raw objectives.
With the global optimization we ascertain that the method will work even if this situation is not given.

\hypertarget{high-precision-operatornameecdfs}{%
\subsection{\texorpdfstring{High-precision \(\operatorname{ECDF}\)s}{High-precision \textbackslash operatorname\{ECDF\}s}}\label{high-precision-operatornameecdfs}}

\begin{figure}

{\centering \includegraphics[width=0.9\linewidth]{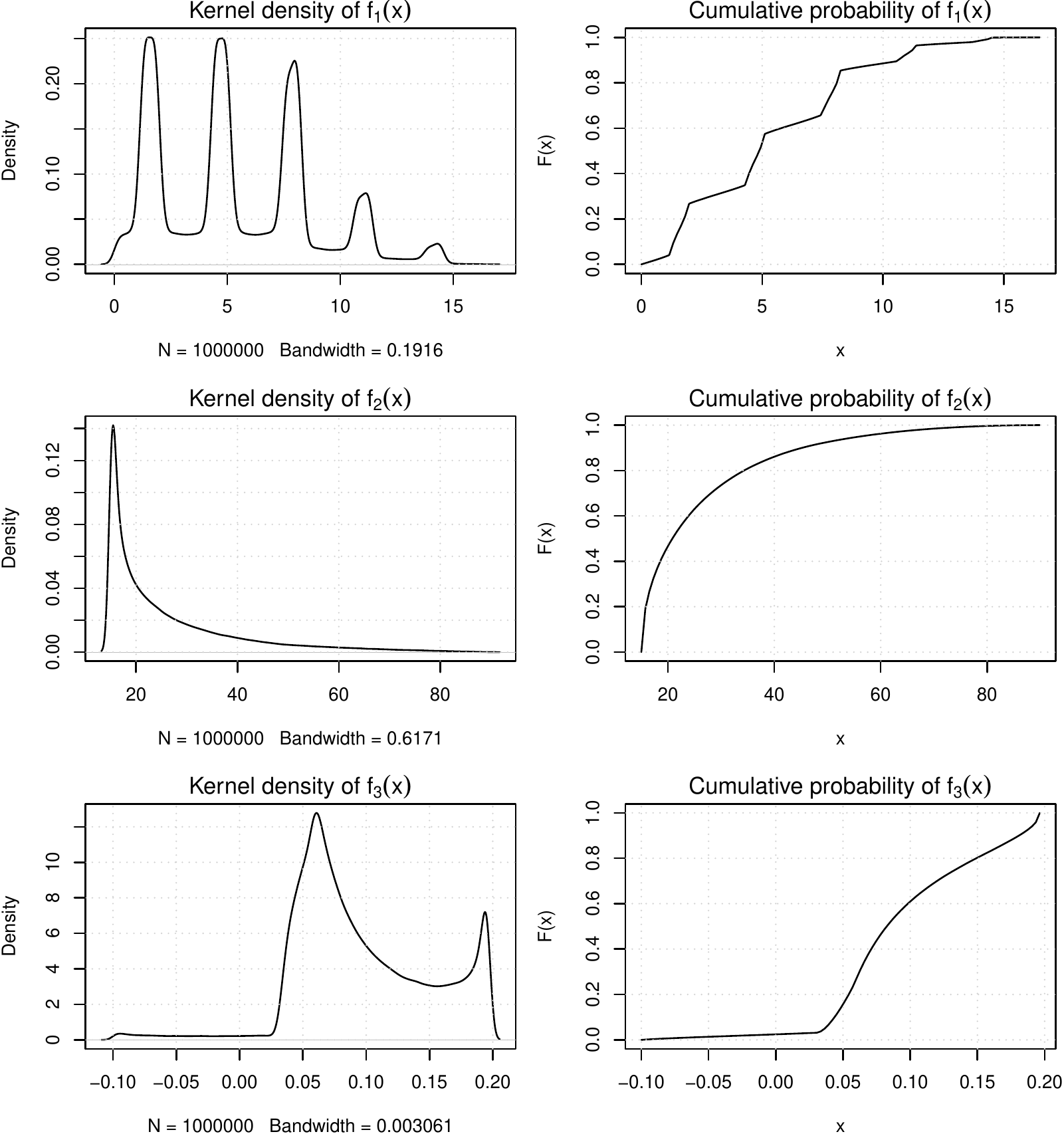} 

}

\caption{The empirical densities and cumulative distributions of the three scaled objectives of the Viennet function.}\label{fig:ex2-high-prec}
\end{figure}

In order to check the validity and precision of the results obtained in this work, we obtain \(10^6\) samples from each objective using random decision vectors.
These are shown in figure \ref{fig:ex2-high-prec}.
Clearly, all of the distributions are non-uniform and different from each other.

\hypertarget{trade-offs}{%
\subsection{Trade-offs}\label{trade-offs}}

Similar to the high-precision \(\operatorname{ECDF}\)s, a priori knowledge about possible trade-offs does usually not exist, except for when the Pareto efficient set has been computed.
Requesting an \emph{infeasible} preference will still lead to an efficient solution, but perhaps not the desired one, either because it does not exist or it cannot be found. This implies that there is some non-linear mapping between the set of infeasible preferences and associated trade-offs.
Ergo, we can split the preference space \(\mathcal{P}\) into the sets of feasible and infeasible trade-offs, denoted by \(\mathcal{P}^\star\) and \(\mathcal{P}^\dagger\), respectively.
If there was access to \(\mathcal{P}^\star\), one could learn about whether some preference \(\mathbf{p}^\star\) leads to the desired trade-off decision vector \(\mathbf{s}\) in the score space. For example, the relationship might be linear.

\hypertarget{results}{%
\section{Results}\label{results}}

In this section, we demonstrate some of the empirical findings.

\hypertarget{ordering-the-solutions}{%
\subsection{Ordering the solutions}\label{ordering-the-solutions}}

In order to introduce order into the set of previously computed Pareto efficient solutions, we first map each solution into the score space and then simply sum up its scores in an unweighted manner, making no additional assumptions about which objective is the most important.

\begin{figure}

{\centering \includegraphics[width=0.9\linewidth]{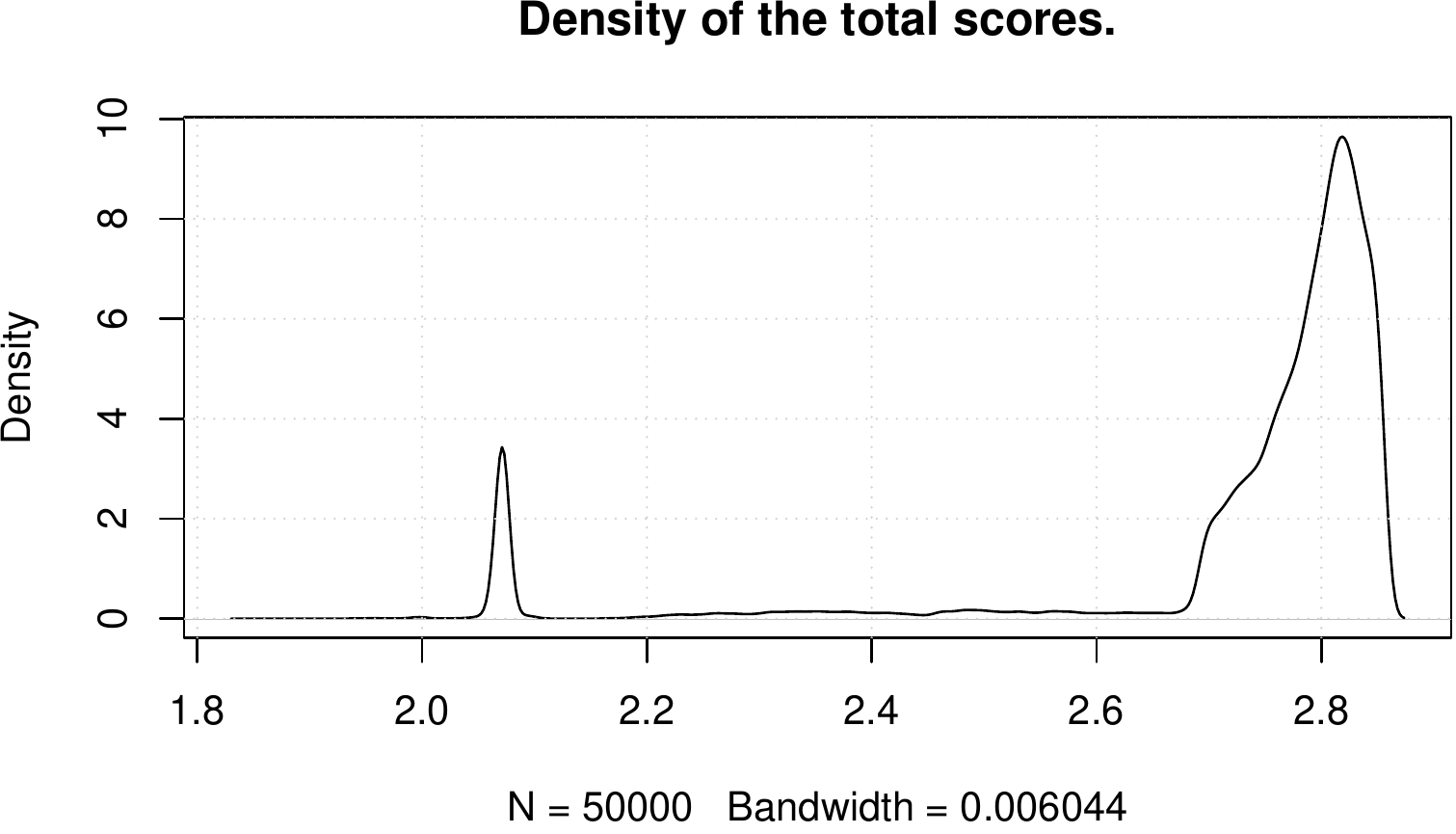} 

}

\caption{The density of the total scores.}\label{fig:dens-total-scores}
\end{figure}

The total scores for each of the solutions are in the interval \(\interval{{\approx}1.848}{{\approx}2.856}\).
This means that any solution with arbitrary combination of scores less than \({\approx}1.848\) is not Pareto efficient.
There is a unique best solution for this problem with objective values of \(\left\{f_1\approx0.9231,f_2\approx15.1532,f_3\approx0.0307\right\}\), marked by a red square in figure \ref{fig:ex2-pareto-front}. The solution marked by the purple square is the worst in the Pareto efficient set.
A lighter color in this figure indicates a higher total score.
The distribution of the total scores is then shown in figure \ref{fig:dens-total-scores}.
In our case, the space containing the highest-scoring solutions is continuous and contiguous, i.e., the Pareto efficient solutions are to be found in a single area, and there is no area disjoint from that which has high scores, too.

\hypertarget{obtaining-desired-trade-offs-1}{%
\subsection{Obtaining desired trade-offs}\label{obtaining-desired-trade-offs-1}}

Obtaining a Pareto efficient solution close to the desired trade-off can be done in multiple ways, all of which require to map the solutions into the score space \(\mathcal{S}\).

\hypertarget{a-priori-correction-of-weights}{%
\subsubsection{A priori correction of weights}\label{a-priori-correction-of-weights}}

\begin{figure}
\centering
\includegraphics{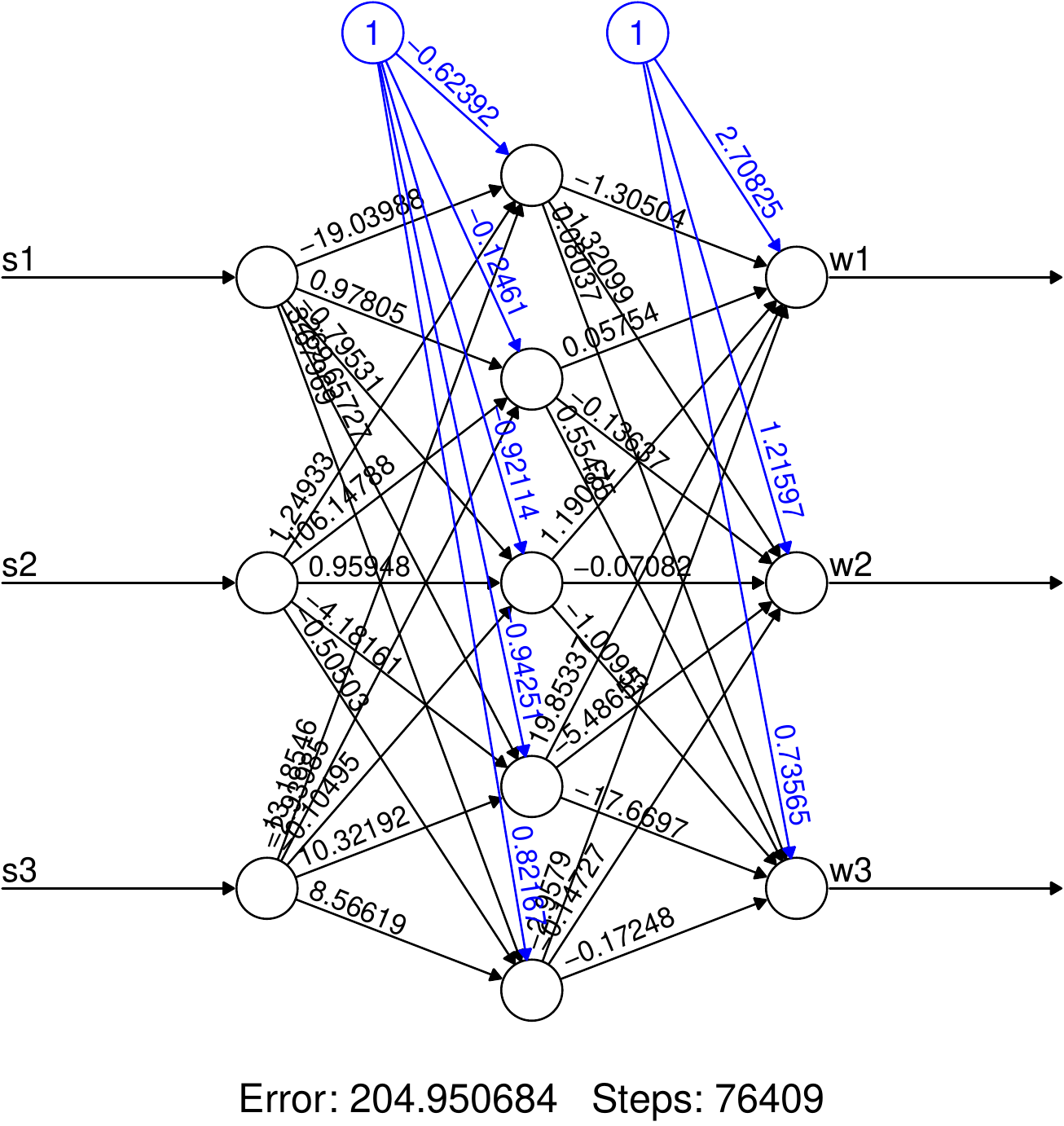}
\caption{\label{fig:ex2-nnet}The neural network for correcting preference.}
\end{figure}

In the first method, we attempt to correct the DMs preference such that it will result in the desired trade-off.
This method requires us to learn a non-linear mapping between the gotten and desired trade-off (\(\mathcal{S}\to\mathcal{P}\)).
This means that some DM expresses their preference in the solution space, which feels more natural. The complexity of choosing appropriate weights for the desired solution is then delegated to this method.
This method does not require to have a notion of whether some expressed preference is feasible. If, however, there is knowledge about \(\mathcal{P}^\star\) or \(\mathcal{P}^\dagger\), and the mapping \(\mathcal{P}^\star\to\mathcal{S}\) is known to be linear, then those preference/solution tuples can be excluded from training.
In order to obtain training data for a model, we require a random subset of the Pareto efficient solutions to learn from.
From the previously approximated Pareto front that holds \(50,000\) efficient solutions, we take now \(2,500\) for this task. The remaining \(47,500\) pairs are held out as validation data and are not in any way part of the training.
To simplify the data, we normalize each solution- and preference-vector such that the most preferred objective has a weight of \(1\).
We attempt to learn an artificial neural network (\protect\hyperlink{ref-rpgk_neuralnet}{Fritsch, Guenther, and Wright 2019}). It is shown in figure \ref{fig:ex2-nnet}.
Since the problem is relatively small, we settle for a single hidden layer with five units. The Sigmoid function is used as activation function.
From that figure we see that the preference is now actually specified in the solution space \(\mathcal{S}\), and that the network will predict the associated preference in the preference space, \(\mathcal{P}\).

\begin{figure}

{\centering \includegraphics[width=0.9\linewidth]{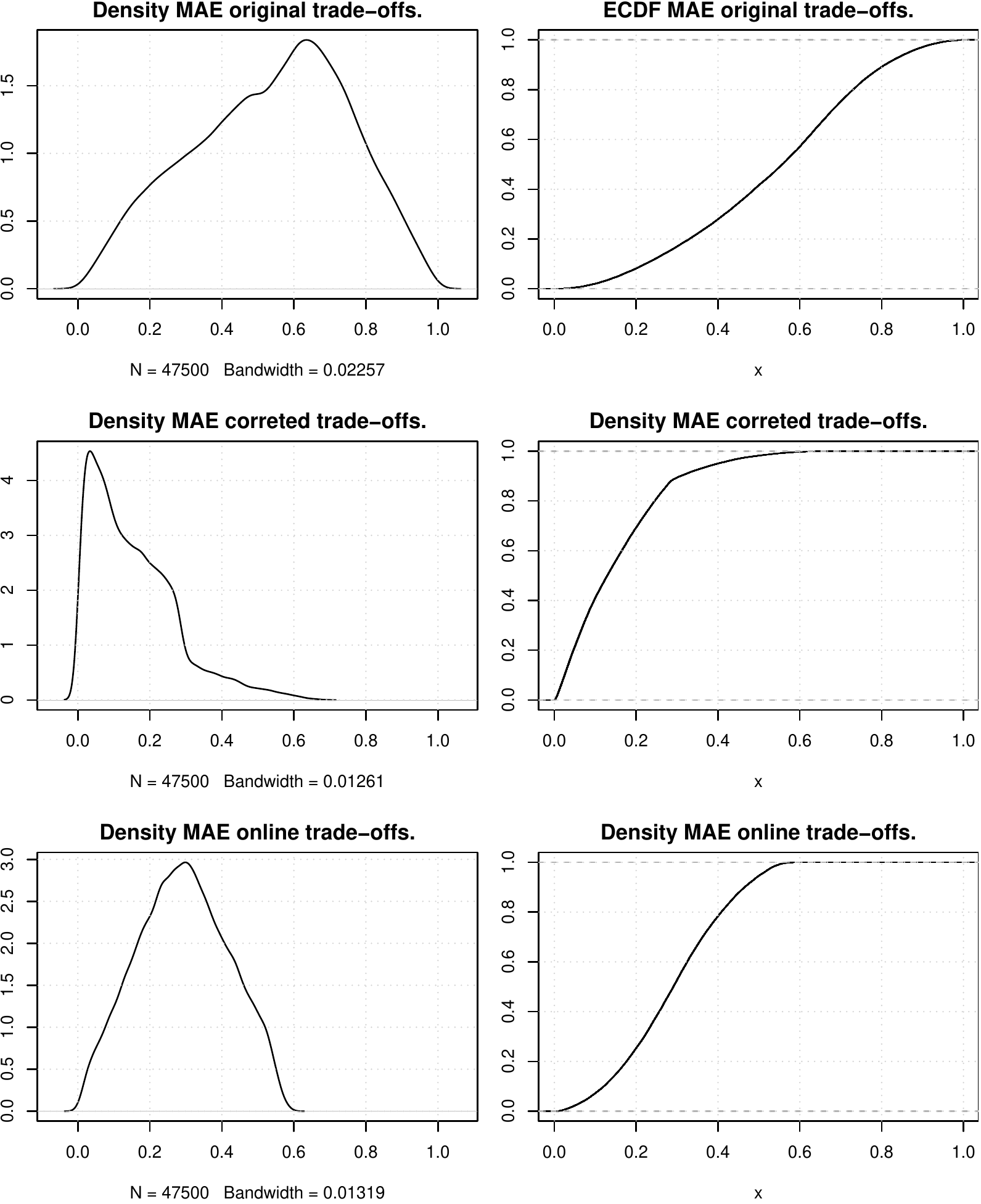} 

}

\caption{The mean absolute error of the original, corrected, and online trade-offs.}\label{fig:ex2-tradeoffs-org-vs-corr-online}
\end{figure}

Figure \ref{fig:ex2-tradeoffs-org-vs-corr-online} shows a significant improvement of the obtained trade-offs using a corrected preference.
We also generate some quantiles (see table \ref{tab:ex2-pref-corr-quant}).
While still not perfect, \(50\)\% of all trade-offs show a mean absolute error of \(\approx0.129\).
We also obtain much less extreme deviations, and no deviation is larger than \(\approx0.679\).

\hypertarget{online-correction}{%
\subsubsection{Online correction}\label{online-correction}}

The method of correcting preference seems promising.
However, we also want to examine the approach of using scores instead of raw objectives during optimization.
If we examine the empirical densities of figure \ref{fig:ex2-high-prec}, it does not appear that we can fit parametric distributions with promising results.
We therefore approximate empirical \(\operatorname{CDF}\)s using \(2,500\) random decision vectors for each objective. While the \(\operatorname{ECDF}\) is a step function, we introduce numerical stability and strict monotonicity by first replacing the step-wise function with a piece-wise linear function. Then, we add slight slopes to the begin and end of each function, to account for extreme values not observed so far. This is conceptually similar to how the utopian vector is created: \(\bm{z}^{\star\star}=\bm{z}^{\star}-\epsilon\).

\begin{table}

\caption{\label{tab:ex2-pref-corr-quant}Quantiles of trade-off errors for the original-, corrected- and online preferences.}
\centering
\begin{tabular}[t]{llll}
\toprule
Quantile & org & corr & online\\
\midrule
0\% & 0.00097 & 0.00014 & 0.0022\\
10\% & 0.22227 & 0.02509 & 0.12085\\
20\% & 0.3299 & 0.0466 & 0.17645\\
30\% & 0.41604 & 0.0702 & 0.21976\\
40\% & 0.48882 & 0.09704 & 0.25631\\
\addlinespace
50\% & 0.55705 & 0.12936 & 0.29087\\
60\% & 0.61541 & 0.16502 & 0.32463\\
70\% & 0.66959 & 0.20302 & 0.36233\\
80\% & 0.72986 & 0.24511 & 0.40791\\
90\% & 0.80762 & 0.30837 & 0.46497\\
\addlinespace
100\% & 0.99669 & 0.67856 & 0.5895\\
\bottomrule
\end{tabular}
\end{table}

Similar to the method of a priori correction of preference, this method achieves trade-offs that are significantly closer to the desired trade-off.
Figure \ref{fig:ex2-tradeoffs-org-vs-corr-online} shows that approx. \(50\)\% or more of all trade-offs are off by less than \(0.3\) for each single objective in the score space.
If we compare the quantiles of this method to the a priori correction method (table \ref{tab:ex2-pref-corr-quant}), then we see that the online method has a consistently larger error than the a prior correction method, but its distribution is more compact, in that the largest error was \(\approx0.09\) less (the \(100\)\% quantile).

\hypertarget{examining-possible-trade-offs}{%
\subsection{Examining possible trade-offs}\label{examining-possible-trade-offs}}

\begin{figure}

{\centering \subfloat[Ideal and homogeneous solution density.\label{fig:ex2-pareto-solution-space-1}]{\includegraphics[width=0.45\linewidth]{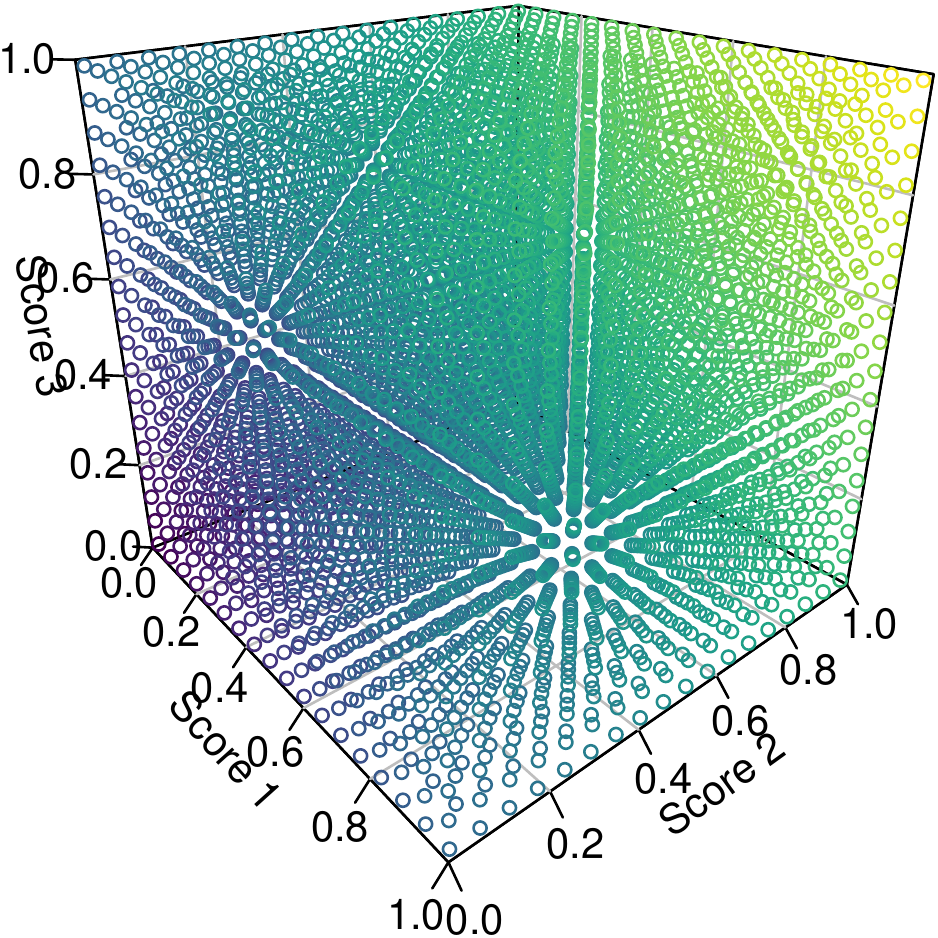} }\subfloat[Inhomogeneous solution density of the Viennet function.\label{fig:ex2-pareto-solution-space-2}]{\includegraphics[width=0.45\linewidth]{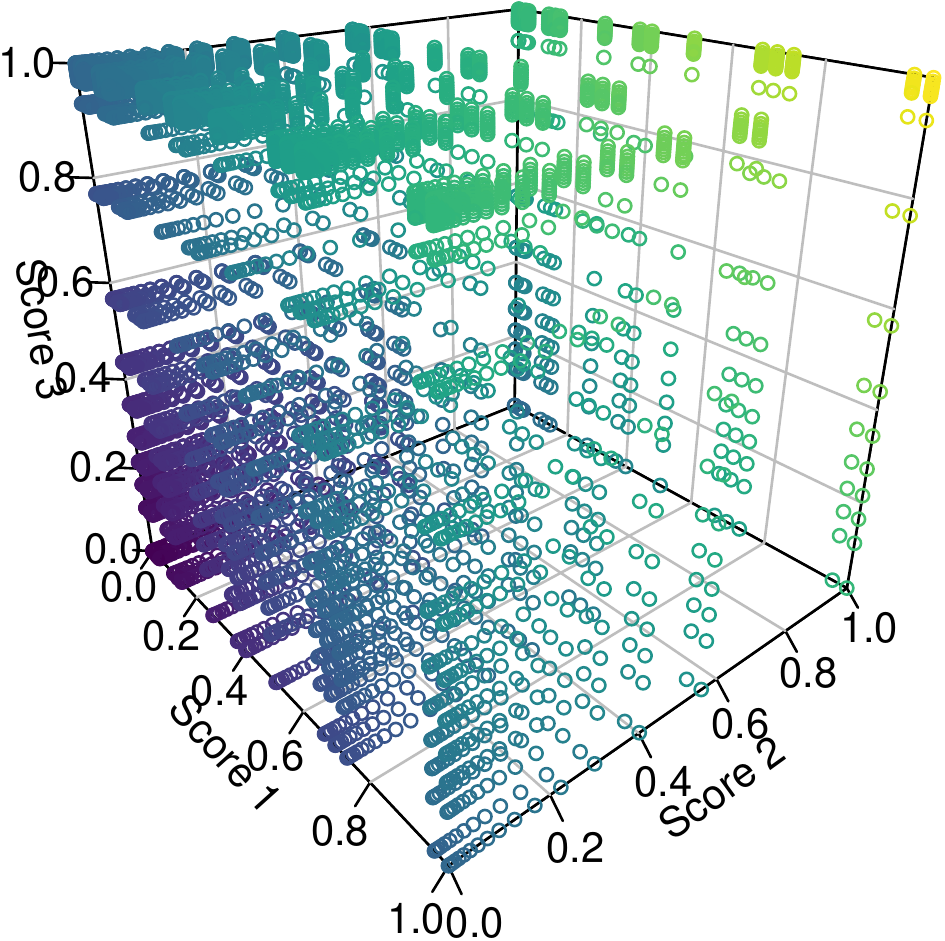} }

}

\caption{Solutions for a three-objective problem.}\label{fig:ex2-pareto-solution-space}
\end{figure}

A DM might expect that the space of all solutions be homogeneous (i.e., having equal density at each point), especially in problems where all objectives have been scaled into a uniform range.
That would imply that each objective's loss is equally likely.
Figure \ref{fig:ex2-pareto-solution-space} (a) demonstrates this assumption, where all points \({\neq}\{1,1,1\}\) represent a linear trade-off between objectives' scores.
This assumption, however, is unlikely to hold in practice. In figure \ref{fig:ex2-pareto-solution-space} (b) we have created a separate linear range for each of the three objectives in the Viennet function, where start and end correspond to lowest and highest observed marginal losses, respectively. Then, these losses were mapped into the score space. Recall that none of the cumulative probability densities from figure \ref{fig:ex2-high-prec} was that of a standard uniform distribution.
Therefore, the mapping results in an inhomogeneous space of feasible (but not necessarily optimal) solutions.
Consider, for example, the second objective. Its marginal losses are in the range of \({\approx}\interval{15}{90}\). The fallacy lies in the expectation that, for example, the two improvements from \(60\) to \(50\) and \(30\) to \(20\) are equally good. However, while the former corresponds to an actual improvement in terms of scores of \({\approx}0.037\), the latter corresponds to a staggering \({\approx}0.269\).
This shows that we can use the mapping into the score space to unveil that a thought-of linearly behaving and homogeneous space of solutions actually is not.
In reality, even low-resolution \(\operatorname{ECDF}\)s can quite significantly help to understand the solution density. Also, their precision improves with each observed sample.

\begin{figure}

{\centering \includegraphics[width=0.9\linewidth]{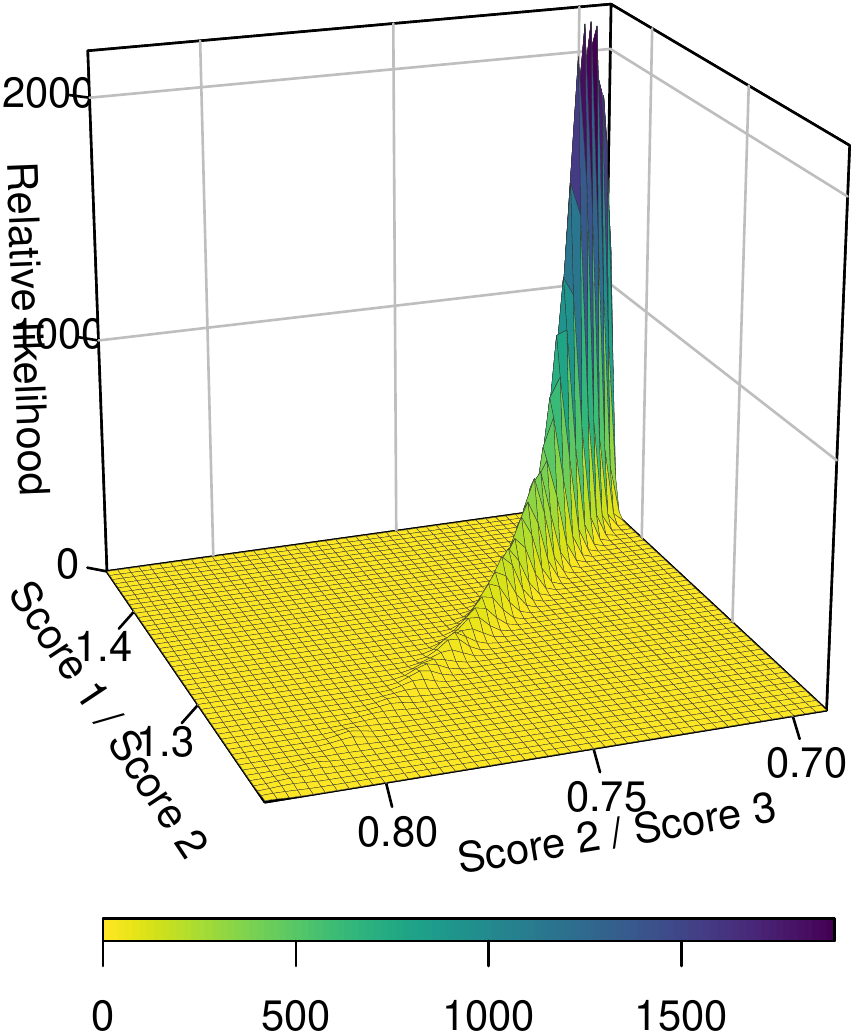} 

}

\caption{The density of the trade-offs that were reached by the solutions in the Pareto efficient set.}\label{fig:ex2-tradeoff-density-2d}
\end{figure}

The Pareto front of the Viennet function has a complex shape. If we map the results gotten into the score space \(\mathcal{S}\), it will allow us to understand which trade-offs are actually feasible.
This is another important piece, as some DM might not be aware that the requested trade-off is not possible (that is, there is no Pareto optimal solution with the given trade-off, or close to it).
Figure \ref{fig:ex2-tradeoff-density-2d} shows the trade-off ratios between the scores that were actually reached by the solutions in the Pareto efficient set.
Without having mapped obtained solutions into the score space, there is no way of understanding whether the desired trade-off was actually obtained.
With the obtained Pareto efficient trade-offs, we can now also attempt to learn whether there is a linear relation between the preference and the associated solution in the score space.

\begin{figure}

{\centering \includegraphics[width=0.9\linewidth]{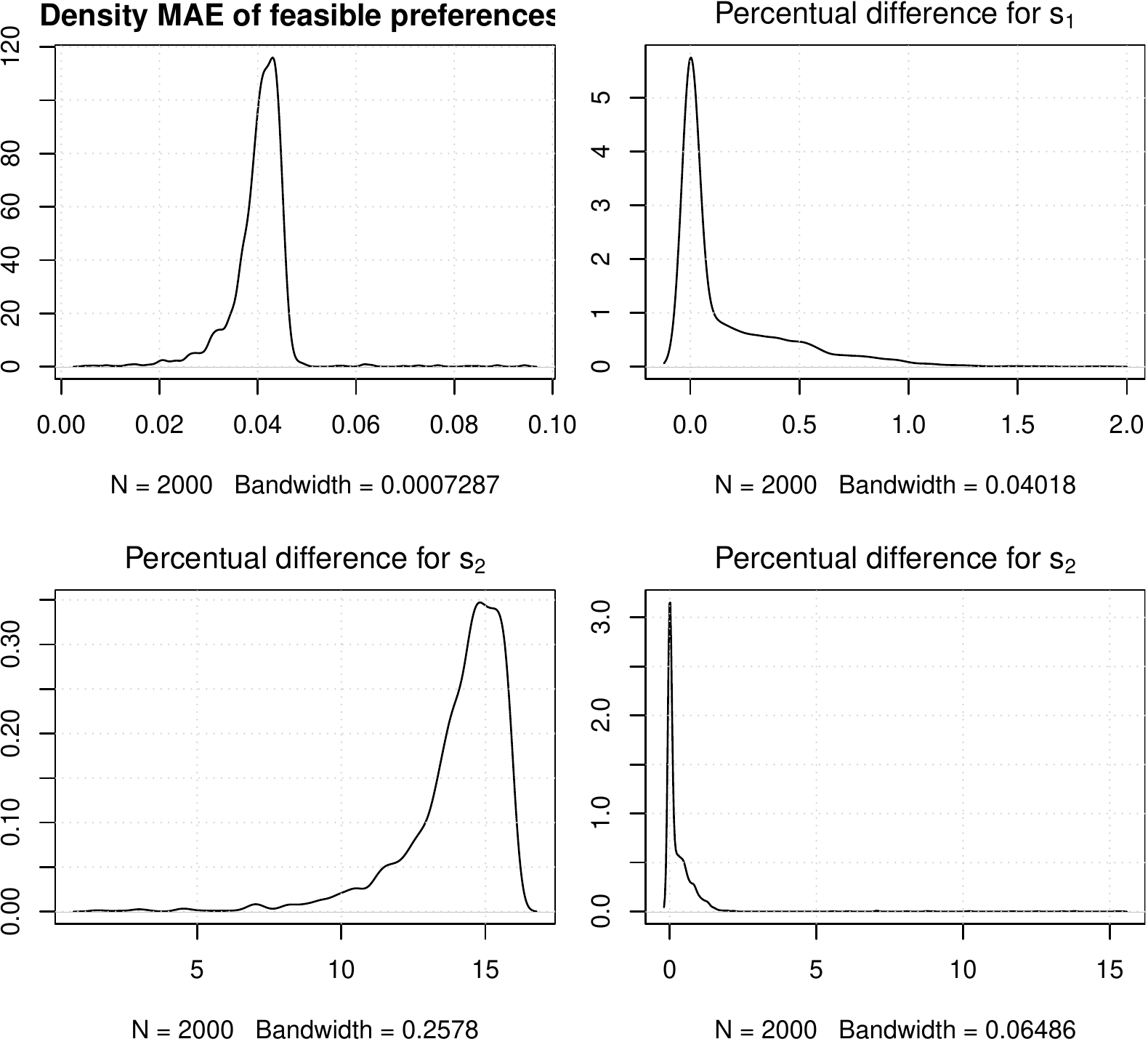} 

}

\caption{Mapping between feasible preferences and solutions in the score space.}\label{fig:ex2-pstar-to-s}
\end{figure}

We have already shown results that indicate that the relation between an expressed preference and the obtained trade-off is of non-linear nature for the Viennet problem.
However, these results included weight preferences for which no or no nearby efficient solutions exist. For example, there are \emph{no} Pareto efficient solutions with a score for the second objective of less than \({\approx}0.681\). Requesting a trade-off that would handicap the second objective to a score less than that is therefore not possible.
Figure \ref{fig:ex2-tradeoff-density-2d} shows the ratios between trade-offs of the first and second score, as well as the second and third score.
We are interested in finding out whether the preferences in the set \(\mathcal{P}^\star\) lead to similar trade-offs, which would imply a linear relationship.
We do not want to reuse any of the solutions generated for the Pareto front, so we will sample inversely from the joint distribution of observed trade-offs.
After running the optimization with the known to be feasible preferences, it turns out that the mapping \(\mathcal{P}^\star\to\mathcal{S}\) is almost linear, as there is only slight deviation.
It also appears that the deviation stems almost exclusively from the second score.
Figure \ref{fig:ex2-pstar-to-s} show the results in more detail.

\hypertarget{references}{%
\section*{References}\label{references}}
\addcontentsline{toc}{section}{References}

\hypertarget{refs}{}
\begin{CSLReferences}{1}{0}
\leavevmode\vadjust pre{\hypertarget{ref-rpgk_neuralnet}{}}%
Fritsch, Stefan, Frauke Guenther, and Marvin N. Wright. 2019. \emph{Neuralnet: Training of Neural Networks}. \url{https://CRAN.R-project.org/package=neuralnet}.

\leavevmode\vadjust pre{\hypertarget{ref-gablonskyK2001}{}}%
Gablonsky, Joerg M., and Carl T. Kelley. 2001. {``A Locally-Biased Form of the {DIRECT} Algorithm.''} \emph{J. Glob. Optim.} 21 (1): 27--37. \url{https://doi.org/10.1023/A:1017930332101}.

\leavevmode\vadjust pre{\hypertarget{ref-miettinen2008}{}}%
Miettinen, Kaisa. 2008. {``Introduction to Multiobjective Optimization: Noninteractive Approaches.''} In \emph{Multiobjective Optimization, Interactive and Evolutionary Approaches}, edited by Jürgen Branke, Kalyanmoy Deb, Kaisa Miettinen, and Roman Slowinski, 5252:1--26. Lecture Notes in Computer Science. Springer. \url{https://doi.org/10.1007/978-3-540-88908-3/_1}.

\leavevmode\vadjust pre{\hypertarget{ref-roy1996theoretical}{}}%
Roy, Bernard, and Vincent Mousseau. 1996. {``A Theoretical Framework for Analysing the Notion of Relative Importance of Criteria.''} \emph{Journal of Multi-Criteria Decision Analysis} 5 (2): 145--59.

\leavevmode\vadjust pre{\hypertarget{ref-viennet1996}{}}%
Viennet, R, Christian Fonteix, and Ivan Marc. 1996. {``Multicriteria Optimization Using a Genetic Algorithm for Determining a Pareto Set.''} \emph{International Journal of Systems Science} 27 (2): 255--60. \url{https://doi.org/10.1080/00207729608929211}.

\leavevmode\vadjust pre{\hypertarget{ref-zeleny1973compromise}{}}%
Zeleny, Milan. 1973. {``Compromise Programming.''} \emph{Multiple Criteria Decision Making}.

\end{CSLReferences}

\end{document}